\documentclass[11pt,a4paper,reqno]{amsart}  
\usepackage{amssymb} \usepackage{amscd} \usepackage{times}
\newcommand{\bea} {\begin{eqnarray}}
\newcommand{\eea} {\end{eqnarray}}
\newcommand{\Bea} {\begin{eqnarray*}}
\newcommand{\Eea} {\end{eqnarray*}}

\def\zbb{\mathbb{Z}}  
  
  \def\phi{\varphi}
 \def\p1{{\mathbb{P}^1_\zbb}}

\newtheorem{Theorem}{\quad Theorem}[section]

\usepackage{amsmath,  amsfonts}

\usepackage[T1]{fontenc}

\newcommand{\be} {\begin{equation}}
\newcommand{\ee} {\end{equation}}

\begin{document}

\title{A compactness result for an equation with Holderian condition.}
\author{Samy Skander Bahoura} 
\address{Departement de Mathematiques, Universite Pierre et Marie Curie, 2 place Jussieu, 75005, Paris, France.}
\email{samybahoura@yahoo.fr, samybahoura@gmail.com}
\maketitle
\begin{abstract}
We give blow-up behavior for a Brezis and Merle's problem with Dirichlet and H\"olderian conditions. Also we derive a compactness creterion as in the work of Brezis and  Merle.\end{abstract}


\section{Introduction and Main Results} 

We set $ \Delta = -(\partial_{11} +\partial_{22}) $  on open analytic domain $ \Omega $ of $ {\mathbb R}^2 $.

\smallskip

We consider the following equation:
\begin{displaymath}  (P)  \left \{ \begin {split} 
      \Delta u  & = V(1+\gamma|x|^{2 \beta}) e^{u}     \,\, &&\text{in} \!\!&&\Omega \subset {\mathbb R}^2, \\
                  u  & = 0  \,\,             && \text{in} \!\!&&\partial \Omega.               
\end {split}\right.
\end{displaymath}

Here, we assume that:

$$  0 \in \partial \Omega, \,\,\, \beta \in [0, 1/2), \,\, \gamma \in [0, \gamma_0], \,\, \gamma_0 >0. $$
and,
$$ 0 \leq V \leq b < + \infty, \,\, e^u \in L^1({\Omega})\,\, {\rm and} \,\,  u \in W_0^{1,1}(\Omega), $$

\smallskip

We can see in [8] a nice formulation of this problem $ (P) $  in the sens of the distributions. This Problem arises from geometrical and physical problems, see for example [1, 3, 21, 24]. The above equation was studied by many authors, with or without  the boundary condition, also for Riemannian surfaces,  see [1-23],  where one can find some existence and compactness results. In [7] we have the following important Theorem,

\smallskip

{\bf Theorem A}{\it (Brezis-Merle [7])}.{\it For $ (u_i)_i $ and $ (V_i)_i $ two sequences of functions relative to $ (P) $ with,
$$ 0 < a \leq V_i \leq b < + \infty $$
then it holds,
$$ \sup_K u_i \leq c, $$
with $ c $ depending on $ a, b, \beta, \gamma_0, K $ and $ \Omega $.}

One can find in [7] an interior estimate if we assume $ a=0 $, but we need an assumption on the integral of $ e^{u_i} $, namely, we have:

\smallskip

{\bf Theorem B}{\it (Brezis-Merle [7])}.{\it For $ (u_i)_i $ and $ (V_i)_i $ two sequences of functions relative to the problem $ (P) $ with,
$$ 0 \leq V_i \leq b < + \infty \,\, {\rm and} \,\, \int_{\Omega} e^{u_i} dy  \leq C, $$
then it holds;
$$ \sup_K u_i \leq c, $$
with $ c $ depending on $ b, \beta, \gamma_0, C, K $ and $ \Omega $.}

\smallskip

We look to the uniform boundedness in all $ \bar \Omega $ of the solutions of the Problem $ (P) $. When $ a=0 $, the boundedness of $ \int_{\Omega} e^{u_i} $ is a necessary condition in the problem $ (P) $ as showed in $ [7] $ by the following counterexample.

\bigskip

{\bf Theorem C}{\it (Brezis-Merle [7])}.{\it There are two sequences $ (u_i)_i $ and $ (V_i)_i $ of the problem $ (P) $ with, 
$$ 0 \leq V_i \leq b < + \infty \,\, {\rm and} \,\, \int_{\Omega} e^{u_i} dy  \leq C, $$
such that,
$$ \sup_{\Omega}  u_i \to + \infty. $$}

To obtain the two first previous results (Theorems A and B) Brezis and Merle used  an inequality (Theorem 1 of [7]) obtained by an approximation argument and used Fatou's lemma and applied the maximum principle in $ W_0^{1,1}(\Omega) $ which arises from Kato's inequality. Also this weak form of the maximum principle is used to prove the local uniform boundedness result by comparing  a certain function and the Newtonian potential. We refer to [6] for a topic about the weak form of the maximum principle.

When $ \gamma = 0 $, the above equation has many properties in the constant and the Lipschitzian cases:

Note that for the problem $ (P) $ ($ \gamma =0 $), by using the Pohozaev identity, we can prove that $ \int_{\Omega} e^{u_i} $ is uniformly bounded when $ 0 < a \leq V_i \leq b < +\infty $ and $  ||\nabla V_i||_{L^{\infty}} \leq A $ and $ \Omega $ starshaped, when $ a=0 $ and $ \nabla \log V_i $ is uniformly bounded, we can bound uniformly $ \int_{\Omega} V_i e^{u_i} $  . In [20], Ma-Wei have proved that those results stay true for all open sets not necessarily starshaped.

\smallskip

In [10] ($ \gamma = 0 $) Chen-Li have proved that if $ a=0 $, $ \nabla \log V_i $ is uniformly bounded and $ u_i $ is locally uniformly bounded in $ L^1 $, then  the functions are uniformly bounded near the boundary.

\smallskip

In [10] ($ \gamma = 0 $) Chen-Li have proved that if $ a=0 $ and $ \int_{\Omega} e^{u_i} $ is uniformly bounded and $ \nabla \log V_i $ is uniformly bounded, then we have the compactness result directly. Ma-Wei in [20], extend this result in the case where $ a >0 $.

\bigskip

If we assume $ V $ more regular, we can have another type of estimates, a $ \sup + \inf $ type inequalities. It was proved by Shafrir see [23], that, if $ (u_i)_i, (V_i)_i $ are two sequences of functions solutions of the previous equation without assumption on the boundary and, $ 0 < a \leq V_i \leq b < + \infty $, then we have the following interior estimate:
$$ C\left (\dfrac{a}{b} \right ) \sup_K u_i + \inf_{\Omega} u_i \leq c=c(a, b, K, \Omega). $$
One can see in [11] an explicit value of $ C\left (\dfrac{a}{b}\right ) =\sqrt {\dfrac{a}{b}} $. In his proof, Shafrir has used a blow-up function, the Stokes formula and an isoperimetric inequality, see [3]. For Chen-Lin, they have used the blow-up analysis combined with some geometric type inequality for the integral curvature.

\bigskip

Now, if we suppose $ (V_i)_i $ uniformly Lipschitzian with $ A $ the
Lipschitz constant, then, $ C(a/b)=1 $ and $ c=c(a, b, A, K, \Omega)
$, see Brezis-Li-Shafrir [5]. This result was extended for
H\"olderian sequences $ (V_i)_i $ by Chen-Lin, see  [11]. Also, one
can see in [18], an extension of the Brezis-Li-Shafrir result to compact
Riemannian surfaces without boundary. One can see in [18] explicit form,
($ 8 \pi m, m\in {\mathbb N}^* $ exactly), for the numbers in front of
the Dirac masses when the solutions blow-up. Here, the notion of isolated blow-up point is used.

In [9] we have some a priori estimates on the 2 and 3-spheres $ {\mathbb S}_2 $, $ {\mathbb S}_3 $.

Here we give the behavior of the blow-up points on the boundary and a proof of Brezis-Merle Problem when $ \gamma \in [0,\gamma_0] $, $ \gamma_0 >0 $ and $ \beta \in [0, 1/2) $.

\smallskip

The Brezis-Merle Problem (see [7]) is:

\smallskip

{\bf Problem}. Suppose that $ V_i \to  V $ in $ C^0( \bar \Omega ) $ with $ 0 \leq V_i \leq b  $ for some positive constant $ b $. Also, we consider a sequence of solutions $ (u_i) $ of $ (P) $ relative to $ (V_i) $ with  $ \gamma =0 $, such that,
$$ \int_{\Omega} e^{u_i} dx \leq C,  $$
is it possible to have:
$$ ||u_i||_{L^{\infty}}\leq C=C(b, C, V, \Omega) ? $$

\bigskip

Here we give a blow-up analysis for a sequence of solutions of the Problem $ (P) $  and a proof of compactness result for the Brezis-Merle's Problem when $ 1/2 > \beta \geq 0 $ and $ \gamma \geq 0 $. We extend the result of Chen-Li [10]. For the blow-up analysis we assume that:
$$ 0 \leq  V_i \leq b, $$
The condition $ V_i \to  V $ in $ C^0(\bar \Omega) $ is not necessary, but for the proof of the compactness result we assume that:
$$ ||\nabla V_i||_{L^{\infty}}\leq  A. $$

\smallskip

Our main result are:

\smallskip

\begin{Theorem}  Assume that $ \max_{\Omega} u_i \to +\infty $, where $ (u_i) $ are solutions of the problem $ (P) $ with:
 $$ \beta \in [0, 1/2),\,\, \gamma \in [0,\gamma_0],\,\, 0 \leq V_i \leq b\,\,\, {\rm and } \,\,\, \int_{\Omega}  e^{u_i} dx \leq C, \,\,\, \forall \,\, i, $$
then, after passing to a subsequence, there is a finction $ u $,  there is a number $ N \in {\mathbb N} $ and $ N  $ points $ x_1, \ldots, x_N \in  \partial \Omega $, such that, 
$$ \partial_{\nu} u_i  \to \partial_{\nu} u +\sum_{j=1}^N \alpha_j \delta_{x_j}, \,\,\, \alpha_j \geq 4\pi, \,\, {\rm weakly\,\, in \, the \, sens \, of \, measures \, on \,\, \partial \Omega.} $$
$$ u_i \to u \,\,\, {\rm in }\,\,\, C^1_{loc}(\bar \Omega-\{x_1,\ldots, x_N \}). $$
 \end{Theorem} 

\smallskip

\begin{Theorem}Assume that $ (u_i) $ are solutions of $ (P) $ relative to $ (V_i) $ with the following conditions:
$$  0 \in \partial \Omega, \,\, \beta \in [0, 1/2), \,\, \gamma \in [0,\gamma_0]. $$
and,
$$ 0 \leq V_i\leq b, \,\,\, ||\nabla V_i||_{L^{\infty}} \leq A \,\,\, {\rm and } \,\,\, \int_{\Omega} e^{u_i} \leq C, $$
we have,
$$  || u_i||_{L^{\infty}} \leq c(b, \beta, \gamma_0, A, C, \Omega), $$
\end{Theorem} 

In the last theorem we extend the result of Chen-Li ($ \gamma =0 $). The proof of Chen-Li and Ma-Wei [10,20], use the moving-plane method ($ \gamma =0 $).

\section{Proof of the theorems} 

\underbar {\it Proof of theorem 1.1:} 

We have:

$$ u_i \in W_0^{1,1}(\Omega). $$

Since $ e^{u_i} \in L^1(\Omega) $ by the corollary 1 of Brezis-Merle's paper (see [7]) we have $ e^{u_i} \in L^k(\Omega) $ for all $ k  >2 $ and the elliptic estimates of Agmon and the Sobolev embedding (see [1]) imply that:

$$ u_i \in W^{2, k}(\Omega)\cap C^{1, \epsilon}(\bar \Omega). $$ 

We denote by $ \partial_{\nu} u_i $ the inner normal derivative. By the maximum principle we have, $ \partial_{\nu} u_i  \geq 0 $.

By the Stokes formula we have, 

$$ \int_{\partial \Omega} \partial_{\nu} u_i d\sigma \leq C, $$

We use the weak convergence in the space of Radon measures to have the existence of a nonnegative Radon measure $ \mu $ such that,

$$ \int_{\partial \Omega} \partial_{\nu} u_i \phi  d\sigma \to \mu(\phi), \,\,\, \forall \,\,\, \phi \in C^0(\partial \Omega). $$

We take an $ x_0 \in \partial \Omega $ such that, $ \mu({x_0}) < 4\pi $. For $ \epsilon >0 $ small enough set $ I_{\epsilon}= B(x_0, \epsilon)\cap \partial \Omega $ on the unt disk or one can assume it as an interval. We choose a function $ \eta_{\epsilon} $ such that,

$$ \begin{cases}
    
\eta_{\epsilon} \equiv 1,\,\,\,  {\rm on } \,\,\,  I_{\epsilon}, \,\,\, 0 < \epsilon < \delta/2,\\

\eta_{\epsilon} \equiv 0,\,\,\, {\rm outside} \,\,\, I_{2\epsilon }, \\

0 \leq \eta_{\epsilon} \leq 1, \\

||\nabla \eta_{\epsilon}||_{L^{\infty}(I_{2\epsilon})} \leq \dfrac{C_0(\Omega, x_0)}{\epsilon}.

\end{cases} $$

We take a $\tilde \eta_{\epsilon} $ such that,

\begin{displaymath}  \left \{ \begin {split} 
      \Delta \tilde \eta_{\epsilon}  & = 0              \,\, &&\text{in} \!\!&&\Omega \subset {\mathbb R}^2, \\
                  \tilde\eta_{\epsilon} & =  \eta_{\epsilon}   \,\,             && \text{in} \!\!&&\partial \Omega.               
\end {split}\right.
\end{displaymath}

{\bf Remark:} We use the following steps in the construction of $ \tilde \eta_{\epsilon} $:

We take a cutoff function $ \eta_{0} $ in $ B(0, 2) $ or $ B(x_0, 2) $:

1- We set $ \eta_{\epsilon}(x)= \eta_0(|x-x_0|/\epsilon) $ in the case of the unit disk it is sufficient.

2- Or, in the general case: we use a chart $ (f, \tilde \Omega) $ with $ f(0)=x_0 $ and we take $ \mu_{\epsilon}(x)= \eta_0 ( f( |x|/ \epsilon)) $ to have  connected  sets $ I_{\epsilon} $ and we take $ \eta_{\epsilon}(y)= \mu_{\epsilon}(f^{-1}(y))$. Because $ f, f^{-1} $ are Lipschitz, $ |f(x)-x_0| \leq k_ 2|x|\leq 1 $ for $ |x| \leq 1/k_2 $ and $ |f(x)-x_0| \geq k_ 1|x|\geq 2 $ for $ |x| \geq 2/k_1>1/k_2 $, the support  of $ \eta $ is in $ I_{(2/k_1)\epsilon} $.

$$ \begin{cases}
    
\eta_{\epsilon} \equiv 1,\,\,\,  {\rm on } \,\,\,  f(I_{(1/k_2)\epsilon}), \,\,\, 0 < \epsilon < \delta/2,\\

\eta_{\epsilon} \equiv 0,\,\,\, {\rm outside} \,\,\, f(I_{(2/k_1)\epsilon }), \\

0 \leq \eta_{\epsilon} \leq 1, \\

||\nabla \eta_{\epsilon}||_{L^{\infty}(I_{(2/k_1)\epsilon})} \leq \dfrac{C_0(\Omega, x_0)}{\epsilon}.

\end{cases} $$

3- Also, we can take: $ \mu_{\epsilon}(x)= \eta_0(|x|/\epsilon) $ and $ \eta_{\epsilon}(y)= \mu_{\epsilon}(f^{-1}(y)) $, we extend it by $ 0 $ outside $ f(B_1(0)) $.  We have $ f(B_1(0)) = D_1(x_0) $, $ f (B_{\epsilon}(0))= D_{\epsilon}(x_0) $ and $ f(B_{\epsilon}^+)= D_{\epsilon}^+(x_0) $ with $ f $ and $ f^{-1} $ smooth diffeomorphism.

$$ \begin{cases}
    
\eta_{\epsilon} \equiv 1,\,\,\,  {\rm on \, a \, the \, connected \, set } \,\,\,  J_{\epsilon} =f(I_{\epsilon}), \,\,\, 0 < \epsilon < \delta/2,\\

\eta_{\epsilon} \equiv 0,\,\,\, {\rm outside} \,\,\, J'_{\epsilon}=f(I_{2\epsilon }), \\

0 \leq \eta_{\epsilon} \leq 1, \\

||\nabla \eta_{\epsilon}||_{L^{\infty}(J'_{\epsilon})} \leq \dfrac{C_0(\Omega, x_0)}{\epsilon}.

\end{cases} $$

And, $ H_1(J'_{\epsilon}) \leq C_1 H_1(I_{2\epsilon}) = C_1 4\epsilon $, since $ f $ is Lipschitz. Here $ H_1 $ is the Hausdorff measure.

We solve the Dirichlet Problem:

\begin{displaymath}  \left \{ \begin {split} 
      \Delta \bar \eta_{\epsilon}  & = \Delta \eta_{\epsilon}              \,\, &&\text{in} \!\!&&\Omega \subset {\mathbb R}^2, \\
                  \bar \eta_{\epsilon} & = 0   \,\,             && \text{in} \!\!&&\partial \Omega.               
\end {split}\right.
\end{displaymath}

and finaly we set $ \tilde \eta_{\epsilon} =-\bar \eta_{\epsilon} + \eta_{\epsilon} $. Also, by the maximum principle and the elliptic estimates we have :

$$ ||\nabla \tilde \eta_{\epsilon}||_{L^{\infty}} \leq C(|| \eta_{\epsilon}||_{L^{\infty}} +||\nabla \eta_{\epsilon}||_{L^{\infty}} + ||\Delta \eta_{\epsilon}||_{L^{\infty}}) \leq \dfrac{C_1}{\epsilon^2}, $$

with $ C_1 $ depends on $ \Omega $.

We use the following estimate, see [4, 8, 14, 25],

$$ ||\nabla u_i||_{L^q} \leq C_q, \,\,\forall \,\, i\,\, {\rm and  }  \,\, 1< q < 2. $$

We deduce from the last estimate that, $ (u_i) $ converge weakly in $ W_0^{1, q}(\Omega) $, almost everywhere to a function $ u \geq 0 $ and $ \int_{\Omega} e^u < + \infty $ (by Fatou lemma). Also, $ V_i $ weakly converge to a nonnegative function $ V $ in $ L^{\infty} $. The function $ u $ is in $ W_0^{1, q}(\Omega) $ solution of :

\begin{displaymath} \left \{ \begin {split} 
      \Delta u  & = V(1+\gamma |x|^{2 \beta}) e^{u} \in L^1(\Omega)    \,\, &&\text{in} \!\!&&\Omega \subset {\mathbb R}^2, \\
                  u  & = 0  \,\,                                     && \text{in} \!\!&&\partial \Omega.               
\end {split}\right.
\end{displaymath}
 
According to the corollary 1 of Brezis-Merle's result, see [7],   we have $ e^{k u }\in L^1(\Omega), k >1 $. By the elliptic estimates, we have $ u \in C^1(\bar \Omega) $.

For two vectors $ f $ and $ g $ we denote by $  f \cdot g $ the inner product of $ f $ and $ g $. 

We can write:

\be \Delta ((u_i-u) \tilde \eta_{\epsilon})= (1+\gamma |x|^{2 \beta})(V_i e^{u_i} -Ve^u)\tilde \eta_{\epsilon} -2\nabla (u_i- u)\cdot \nabla \tilde \eta_{\epsilon}. \label{(1)}\ee

We use the interior esimate of Brezis-Merle, see [7],

\bigskip

\underbar {\it Step 1:} Estimate of the integral of the first term of the right hand side of $ (\ref{(1)}) $.

\bigskip

We use the Green formula between $ \tilde \eta_{\epsilon} $ and $ u $, we obtain,

\be  \int_{\Omega} (1+\gamma |x|^{2 \beta})Ve^u \tilde \eta_{\epsilon} dx =\int_{\partial \Omega} \partial_{\nu} u \eta_{\epsilon} \leq C'\epsilon ||\partial_{\nu}u||_{L^{\infty}}= C \epsilon \label{(2)}\ee

We have,

\begin{displaymath} \left \{ \begin {split} 
      \Delta u_i  & = (1+\gamma |x|^{2 \beta})V_i e^{u_i}                      \,\, &&\text{in} \!\!&&\Omega \subset {\mathbb R}^2, \\
                  u_i  & = 0  \,\,                                     && \text{in} \!\!&&\partial \Omega.               
\end {split}\right.
\end{displaymath}

We use the Green formula between $ u_i $ and $ \tilde \eta_{\epsilon} $ to have:

\be \int_{\Omega} (1+\gamma |x|^{2 \beta})V_i e^{u_i} \tilde \eta_{\epsilon} dx = \int_{\partial \Omega} \partial_{\nu} u_i \eta_{\epsilon} d\sigma \to \mu(\eta_{\epsilon}) \leq \mu(J'_{\epsilon}) \leq 4  \pi - \epsilon_0, \,\,\, \epsilon_0 >0 \label{(3)}\ee

From $ (\ref{(2)}) $ and $ (\ref{(3)}) $ we have for all $ \epsilon >0 $ there is $ i_0 =i_0(\epsilon) $ such that, for $ i \geq i_0 $,

\be \int_{\Omega} |(1+\gamma |x|^{2 \beta})(V_ie^{u_i}-Ve^u) \tilde \eta_{\epsilon}| dx \leq 4 \pi -\epsilon_0 + C \epsilon \label{(4)}\ee

\bigskip

\underbar {\it Step 2:} Estimate of integral of the second term of the right hand side of $ (\ref{(1)}) $.

\bigskip

Let $ \Sigma_{\epsilon} = \{ x \in \Omega, d(x, \partial \Omega) = \epsilon^3 \} $ and $ \Omega_{\epsilon^3} = \{ x \in \Omega, d(x, \partial \Omega) \geq \epsilon^3 \} $, $ \epsilon > 0 $. Then, for $ \epsilon $ small enough, $ \Sigma_{\epsilon} $ is hypersurface.

The measure of $ \Omega-\Omega_{\epsilon^3} $ is $ k_2\epsilon^3 \leq meas(\Omega-\Omega_{\epsilon^3})= \mu_L (\Omega-\Omega_{\epsilon^3}) \leq k_1 \epsilon^3 $.

{\bf Remark}: for the unit ball $ {\bar B(0,1)} $, our new manifold is $ {\bar B(0, 1-\epsilon^3)} $.

( Proof of this fact; let's consider $ d(x, \partial \Omega) = d(x, z_0), z_0 \in \partial \Omega $, this imply that $ (d(x,z_0))^2 \leq (d(x, z))^2 $ for all $ z \in \partial \Omega $ which it is equivalent to $ (z-z_0) \cdot (2x-z-z_0) \leq 0 $ for all $ z \in \partial \Omega $, let's consider a chart around $ z_0 $ and $ \gamma (t) $ a curve in $ \partial \Omega $, we have;

$ (\gamma (t)-\gamma(t_0) \cdot (2x-\gamma(t)-\gamma(t_0)) \leq 0 $ if we divide by $ (t-t_0) $ (with the sign and tend $ t $ to $ t_0 $), we have $ \gamma'(t_0) \cdot (x-\gamma(t_0)) = 0 $, this imply that $ x= z_0-s \nu_0 $ where $ \nu_0 $ is the outward normal of $ \partial \Omega $ at $ z_0 $))

With this fact, we can say that $ S= \{ x, d(x, \partial \Omega) \leq \epsilon \}= \{ x= z_0- s \nu_{z_0}, z_0 \in \partial \Omega, \,\, -\epsilon \leq s \leq \epsilon \} $. It  is sufficient to work on  $ \partial \Omega $. Let's consider a charts $ (z, D=B(z, 4 \epsilon_z), \gamma_z) $ with $ z \in \partial \Omega $ such that $ \cup_z B(z, \epsilon_z) $ is  cover of $ \partial \Omega $ .  One can extract a finite cover $ (B(z_k, \epsilon_k)), k =1, ..., m $, by the area formula the measure of $ S \cap B(z_k, \epsilon_k) $ is less than a $ k\epsilon $ (a $ \epsilon $-rectangle).  For the reverse inequality, it is sufficient to consider one chart around one point of the boundary.

We write,

\be \int_{\Omega} |\nabla ( u_i -u) \cdot \nabla \tilde \eta_{\epsilon} | dx =
\int_{\Omega_{\epsilon^3}} |\nabla (u_i -u) \cdot \nabla \tilde \eta_{\epsilon}| dx + \int_{\Omega - \Omega_{\epsilon^3}} |\nabla (u_i-u) \cdot \nabla \tilde \eta_{\epsilon}| dx.  \label{(5)}\ee

\bigskip

\underbar {\it Step 2.1:} Estimate of $ \int_{\Omega - \Omega_{\epsilon^3}} |\nabla (u_i-u) \cdot \nabla \tilde \eta_{\epsilon}| dx $.

\bigskip

First, we know from the elliptic estimates that  $ ||\nabla \tilde \eta_{\epsilon}||_{L^{\infty}} \leq C_1 /\epsilon^2 $, $ C_1 $ depends on $ \Omega $

We know that $ (|\nabla u_i|)_i $ is bounded in $ L^q, 1< q < 2 $, we can extract  from this sequence a subsequence which converge weakly to $ h \in L^q $. But, we know that we have locally the uniform convergence to $ |\nabla u| $ (by Brezis-Merle's theorem), then, $ h= |\nabla u| $ a.e. Let $ q' $ be the conjugate of $ q $.

We have, $  \forall f \in L^{q'}(\Omega)$

$$ \int_{\Omega} |\nabla u_i| f dx \to \int_{\Omega} |\nabla u| f dx $$

If we take $ f= 1_{\Omega-\Omega_{\epsilon^3}} $, we have:

$$ {\rm for } \,\, \epsilon >0 \,\, \exists \,\, i_1 = i_1(\epsilon) \in {\mathbb N}, \,\,\, i \geq  i_1,  \,\, \int_{\Omega-\Omega_{\epsilon^3}} |\nabla u_i| \leq \int_{\Omega-\Omega_{\epsilon^3}} |\nabla u| + \epsilon^3. $$

Then, for $ i \geq i_1(\epsilon) $,

$$ \int_{\Omega-\Omega_{\epsilon^3}} |\nabla u_i| \leq meas(\Omega-\Omega_{\epsilon^3}) ||\nabla u||_{L^{\infty}} + \epsilon^3 = \epsilon^3(k_1 ||\nabla u||_{L^{\infty}} + 1). $$

Thus, we obtain,

\be \int_{\Omega - \Omega_{\epsilon^3}} |\nabla (u_i-u) \cdot \nabla \tilde \eta_{\epsilon}| dx \leq  \epsilon C_1(2 k_1 ||\nabla u||_{L^{\infty}} + 1) \label{(6)}\ee

The constant $ C_1 $ does  not depend on $ \epsilon $ but on $ \Omega $.

\bigskip

\underbar {\it Step 2.2:} Estimate of $ \int_{\Omega_{\epsilon^3}} |\nabla (u_i-u) \cdot \nabla \tilde \eta_{\epsilon}| dx $.

\bigskip

We know that, $ \Omega_{\epsilon} \subset \subset \Omega $, and ( because of Brezis-Merle's interior estimates) $ u_i \to u $ in $ C^1(\Omega_{\epsilon^3}) $. We have,

$$ ||\nabla (u_i-u)||_{L^{\infty}(\Omega_{\epsilon^3})} \leq \epsilon^3,\, {\rm for } \,\, i \geq i_3 = i_3(\epsilon). $$

We write,
 
$$ \int_{\Omega_{\epsilon^3}} |\nabla (u_i-u) \cdot \nabla \tilde \eta_{\epsilon}| dx \leq ||\nabla (u_i-u)||_{L^{\infty}(\Omega_{\epsilon^3})} ||\nabla \tilde \eta_{ \epsilon}||_{L^{\infty}} \leq C_1 \epsilon \,\, {\rm for } \,\, i \geq i_3, $$

For $ \epsilon >0 $, we have for $ i \in {\mathbb N} $, $ i \geq \max \{i_1, i_2, i_3 \} $,

\be \int_{\Omega} |\nabla (u_i-u) \cdot \nabla \tilde \eta_{\epsilon}| dx \leq \epsilon C_1(2 k_1 ||\nabla u||_{L^{\infty}} + 2) \label{(7)}\ee

From $ (\ref{(4)}) $ and $ (\ref{(7)}) $, we have, for $ \epsilon >0 $, there is $ i_3= i_3(\epsilon) \in {\mathbb N}, i_3 = \max \{ i_0, i_1, i_2 \} $ such that,

\be \int_{\Omega} |\Delta [(u_i-u)\tilde \eta_{\epsilon}]|dx \leq 4 \pi-\epsilon_0+  \epsilon 2 C_1(2 k_1 ||\nabla u||_{L^{\infty}} + 2 + C) \label{(8)}\ee

We choose $ \epsilon >0 $ small enough to have a good estimate of  $ (\ref{(1)}) $.

Indeed, we have:

\begin{displaymath} \left \{ \begin {split} 
      \Delta [(u_i-u) \tilde \eta_{\epsilon}]   & = g_{i,\epsilon}                   \,\, &&\text{in} \!\!&&\Omega \subset {\mathbb R}^2, \\
                 (u_i-u) \tilde \eta_{\epsilon}    & = 0  \,\,                                     && \text{in} \!\!&&\partial \Omega.               
\end {split}\right.
\end{displaymath}

with $ ||g_{i, \epsilon} ||_{L^1(\Omega)} \leq 4 \pi -\dfrac{\epsilon_0}{2}. $

We can use Theorem 1 of [7] to conclude that there are $ q \geq \tilde q >1 $ such that:

$$ \int_{V_{\epsilon}(x_0)} e^{\tilde q |u_i-u|} dx \leq \int_{\Omega} e^{q|u_i -u| \tilde \eta_{\epsilon}} dx \leq C(\epsilon,\Omega). $$
 
where, $ V_{\epsilon}(x_0) $ is a neighberhood of $ x_0 $ in $ \bar \Omega $. Here we have used that in a neighborhood of $ x_0 $  by the elliptic estimates, 
$ 1- C \epsilon \leq \tilde \eta_{\epsilon} \leq 1 $. (We can take $ B(x_0,\epsilon^3) $).

Thus, for each $ x_0 \in \partial \Omega - \{ \bar x_1,\ldots, \bar x_m \} $ there is $ \epsilon_{x_0} >0, q_{x_0} > 1 $ such that:

\be \int_{B(x_0, \epsilon_{x_0})} e^{q_{x_0} u_i} dx \leq C, \,\,\, \forall \,\,\, i. \label{(9)}\ee

Now, we consider a cutoff function $ \eta \in C^{\infty}({\mathbb R}^2) $ such that

$$ \eta \equiv 1 \,\,\, {\rm on } \,\,\, B(x_0, \epsilon_{x_0}/2) \,\,\, {\rm and } \,\,\, \eta \equiv 0 \,\,\, {\rm on } \,\,\, {\mathbb R}^2 -B(x_0, 2\epsilon_{x_0}/3). $$

We write

$$ -\Delta (u_i \eta) = (1+\gamma |x|^{2\beta})V_i e^{u_i} \eta - 2 \nabla u_i \cdot \nabla \eta  - u_i \Delta \eta. $$

By the elliptic estimates (see [15]) $ (u_i)_i $ is uniformly bounded in $ W^{2, q_1}(V_{\epsilon}(x_0)) $ and also, in $ C^1(V_{\epsilon}(x_0)) $. Finaly, we have, for some $ \epsilon > 0 $ small enough,

$$ || u_i||_{C^{1,\theta}[B(x_0, \epsilon)]} \leq c_3 \,\,\, \forall \,\,\, i. $$

We have proved that, there is a finite number of points $ \bar x_1, \ldots, \bar x_m $ such that the squence $ (u_i)_i  $ is locally uniformly bounded (in $ C^{1,\theta}, \theta >0 $) in $ \bar \Omega - \{ \bar x_1, \ldots , \bar x_m \} $.

\bigskip

\underbar {\it Proof of theorem 1.2:} 

\bigskip

Without loss of generality, we can assume that $ 0 $ is a blow-up point. Since the boundary is an analytic curve $ \gamma_1(t) $, there is a neighborhood of  $ 0 $ such that the curve $ \gamma_1 $ can be extend to a holomorphic map such that $ \gamma_1'(0) \not = 0 $ (series) and by the inverse mapping one can assume that this map is univalent around $ 0 $. In the case when the boundary is a simple Jordan curve the domain is simply connected, see [24]. In the case that the domains has a finite number of holes it is conformally equivalent to a disk with a finite number of disks removed, see [17]. Here we consider a general domain. Without loss of generality one can assume that $ \gamma_1 (B_1^+) \subset \Omega $ and also $ \gamma_1 (B_1^-) \subset (\bar \Omega)^c $ and $ \gamma_1 (-1,1) \subset \partial \Omega $ and $ \gamma_1 $ is univalent. This means that $ (B_1, \gamma_1) $ is a local chart around $ 0 $ for $ \Omega $ and $ \gamma_1 $ univalent. (This fact holds if we assume that we have an analytic domain, in the sense of Hofmann see [16], (below a graph of an analytic function), we have necessary the condition $ \partial \bar \Omega = \partial \Omega $ and the graph is analytic, in this case $ \gamma_1 (t)= (t, \phi(t)) $ with $ \phi $ real analytic and an example of this fact is the unit disk  around the point $ (0,1) $ for example).

By this conformal transformation, we can assume that $ \Omega =B_1^+ $, the half ball, and $ \partial^+ B_1^+ $ is the exterior part, a part which not contain $ 0 $ and on which  $ u_i $ converge in the $ C^1 $ norm to $ u $. Let us consider $ B_{\epsilon}^+ $, the half ball with radius $ \epsilon >0 $. Also, one can consider a $ C^1 $ domain (a rectangle between two half disks) and by charts its image is a $ C^1 $ domain)
We know that:

$$ u_i \in W^{2, k}(\Omega)\cap C^{1, \epsilon}(\bar \Omega). $$ 

Thus we can use integrations by parts (Stokes formula). The second Pohozaev identity applied around the blow-up $ 0 $ see for example [2, 20, 22] gives :

\be  \int_{B_{\epsilon}^+} \Delta u_i (x \cdot \nabla u_i) dx = -\int_{\partial^+ B_{\epsilon}^+}  g(\nabla u_i)d\sigma, \label{(10)}\ee

with,

$$ g(\nabla u_i)= (\nu \cdot \nabla u_i)(x \cdot\nabla u_i)- x \cdot \nu \dfrac{|\nabla u_i|^2}{2}. $$

Thus,

\be  \int_{B_{\epsilon}^+} V_i (1+\gamma |x|^{2 \beta}) e^{u_i} (x \cdot \nabla u_i) dx =-\int_{\partial^+ B_{\epsilon}^+}  g(\nabla u_i)d\sigma, \label{(10)}\ee

After integration by parts, we obtain:

$$  \int_{B_{\epsilon}^+} 2V_i(1+ (1+ \beta)\gamma |x|^{2 \beta}) e^{u_i} dx +  \int_{B_{\epsilon}^+} x \cdot\nabla V_i (1+\gamma |x|^{2 \beta})e^{u_i} dx - \int_{\partial B_{\epsilon}^+} \nu \cdot x (1+\gamma |x|^{2 \beta})V_ie^{u_i} d\sigma = $$

\be =  \int_{\partial^+ B_{\epsilon}^+} g(\nabla u_i)d\sigma, \label{(10)} \ee

Also, for $ u $  we have:

$$  \int_{B_{\epsilon}^+} 2V(1+(1+ \beta)\gamma |x|^{2 \beta})e^{u} dx +  \int_{B_{\epsilon}^+} x \cdot \nabla V (1+\gamma |x|^{2 \beta})e^{u} dx - \int_{\partial B_{\epsilon}^+}  \nu \cdot x(1+ |x|^{2 \beta})Ve^{u}d\sigma = $$

\be = \int_{\partial^+ B_{\epsilon}^+} g(\nabla u)d\sigma, \ee

We use the fact that $ u_i=u=0 $ on $ \{x_1=0 \} $ and $ u_i, u $ are bounded in the $ C^1 $ norm outside a neighborhood of $ 0 $ and we tend $ i $ to $ +\infty $ and then $ \epsilon $ to  $ 0 $ to obtain:

\be \int_{B_{\epsilon}^+} V_i (1+\gamma |x|^{2 \beta})e^{u_i} dx = o(1)+O(\epsilon), \ee

however

\be \int_{ \gamma_1 ( B_{\epsilon}^+)} V_i (1+\gamma |x|^{2 \beta})e^{u_i} dx =\int_{\partial  \gamma_1 ( B_{\epsilon}^+)} \partial_{\nu} u_i d\sigma = \alpha_1+ O(\epsilon)+o(1) >0. \ee

which is a contradiction.

Here we used a theorem of Hofmann see [16], which gives the fact that $ \gamma_1 (B_{\epsilon}^+) $ is a Lipschitz domain. Also, we can see that $ \gamma_1 ((-\epsilon, \epsilon)) $ and $ \gamma_1 (\partial^+ B_{\epsilon}^+) $ are submanifolds.  

We start with a Lipschitz domain $ B_{\epsilon}^+ $ because it is convex and by the univalent and conformal map $ \gamma $ the image of this domain $ \gamma_1 (B_{\epsilon}^+) $ is a Lipschitz domain and thus we can apply the integration by part and here we know the explicit formula of the unit outward normal it is the usual unit outward normal (normal to the tangent space of the boundary which we know explicitly because we have two submanifolds).

In the case of the disk $ D = \Omega $, it is sufficient to consider $ B(0,\epsilon)  \cap D $ which is a Lipschitz domain because it is convex (and not necessarily $ \gamma_1 (B_{\epsilon}^+) $).

There is a version of the integration by part which is the Green-Riemann formula in dimension 2 on a domain $ \Omega $. This formula holds if we assume that there is a finite number of points $ y_1,..., y_m $ such that $ \partial \Omega - (y_1,..., y_m) $ is a $ C^1 $ manifold, see [2], for the Gauss-Green-Riemann-Stokes formula, for $ C^1 $ domains with singular points (here a finite number of singular points).

Remark: Note that a monograph of Droniou contain a proof of all fact about Sobolev spaces (with Strong Lipschitz property) with only weak Lipschitz property (Lipschitz-Charts), we start with Strong Lipschitz property and by $ \gamma_1 $ we have weak Lipschtz property.

\end{document}